\documentclass[11pt]{amsart}

\usepackage[dvips]{graphicx}
\usepackage{color}              
\usepackage{epsfig}

\newtheorem{theorem}{Theorem}[section]
\newtheorem{proposition}[theorem]{Proposition}
\newtheorem{corollary}[theorem]{Corollary}
\newtheorem{lemma}[theorem]{Lemma}

\theoremstyle{definition}

\theoremstyle{remark}
\newtheorem{remark}[theorem]{Remark}

\newcommand{\thmref}[1]{Theorem~\ref{#1}}
\newcommand{\propref}[1]{Proposition~\ref{#1}}
\newcommand{\secref}[1]{\S\ref{#1}}
\newcommand{\lemref}[1]{Lemma~\ref{#1}}
\newcommand{\corref}[1]{Corollary~\ref{#1}} 
\newcommand{\figref}[1]{Fig.~\ref{#1}}

\newcommand{\remref}[1]{Remark~\ref{#1}}

\DeclareMathOperator{\short}{short}
\DeclareMathOperator{\thin}{Thin}
\DeclareMathOperator{\Mod}{Mod}
\DeclareMathOperator{\Ext}{Ext}
\DeclareMathOperator{\area}{area}

\newcommand{\co}{\colon\thinspace}
\newcommand{\emul}{\overset{.}{\asymp}}
\newcommand{\gmul}{\overset{.}{\succ}}
\newcommand{\lmul}{\overset{.}{\prec}}
\newcommand{\eadd}{\overset{+}{\asymp}}
\newcommand{\gadd}{\overset{+}{\succ}}
\newcommand{\ladd}{\overset{+}{\prec}}

\newcommand{\CC}{{\mathcal C}}

\newcommand{\T}{{\mathcal T}}
\newcommand{\M}{{\mathcal M}}

\newcommand{\Teich}{Teichm\"uller }
\newcommand{\hyp}{{\mathbb H}}

\newcommand\Z{{\mathbb Z}}
\newcommand\R{{\mathbb R}}

\newcommand{\ep}{\epsilon}
\newcommand{\npm}{\nu_{\pm}}
\newcommand{\np}{\nu_+}
\newcommand{\nm}{\nu_-}

\newcommand{\npb}{\bar \nu_+}
\newcommand{\nmb}{\bar \nu_-}

\newcommand{\p}{\partial}

\begin{document}

\title{A combinatorial model for the Teichm{\"u}ller metric}
\author{Kasra Rafi}
\date{\today}
\maketitle

\section{Introduction} \label{sec:intro}
\noindent
This paper should be considered a sequel to \cite{me:SC}. We continue here 
to study the geometry of \Teich space using combinatorial properties of curves 
on surfaces. The main result is a formula for the \Teich distance between two 
points in \Teich space, in terms of the combinatorial information extracted from 
short curves of these two points. Let $S$ be a surface of finite type with negative 
Euler characteristic and let $\sigma_1$ and $\sigma_2$ be two points in the thick 
part of \Teich space $\T(S)$ of $S$. Let $\mu_1$ and $\mu_2$ be short markings 
on $\sigma_1$ and $\sigma_2$, respectively. 

\begin{theorem} \label{thm:main}
There exists $k>0$ such that 
\begin{equation}
d_\T(\sigma_1, \sigma_2 ) \asymp \sum_Y \big[ d_Y(\mu_1, \mu_2)\big]_k +
\sum_\alpha \log \big[ d_\alpha(\mu_1,\mu_2) \big]_k.
\end{equation}
\end{theorem}

In the above theorem, the first sum is over all subsurfaces of $S$ that are not 
annuli and the second sum is over all simple closed curves on $S$;
$d_Y(\mu_1,\mu_2)$ measures the relative complexity of the restrictions
of $\mu_1$ and $\mu_2$ to a subsurface $Y$, and $d_\alpha(\mu_1,\mu_2)$
measures the relative twisting of $\mu_1$ and $\mu_2$ around a curve 
$\alpha$; the function $[x]_k$ is equal to zero when $x < k$ and is equal to $x$ 
when $x \geq k$, that is, we take into account only terms that are large enough;
and the function $\log$ is a modified logarithm so that, for $x \in [0,1]$, 
$\log x = 0$. A general version of this theorem, where $\sigma_1$ and 
$\sigma_2$ are not necessarily in the thick part, is stated in \secref{sec:general}
(\thmref{thm:general}).

Other recent results relate the geometry of \Teich space to combinatorial 
spaces. In \cite{minsky:CCI} Masur and Minsky show that the electrified 
\Teich space is quasi-isometric to the complex of curves and therefore is 
also $\delta$--hyperbolic.
Brock has shown (\cite{brock:WV}) that \Teich space equipped with the 
Weil-Petersson metric is quasi-isometric to the pants complex. 
Most recent developments in 
studying the Weil-Petersson metric have resulted from this analogy.

To drive our formula, we need to acquire an understanding of 
how the length and the twisting parameter of a curve change along a 
\Teich geodesic. \cite{me:SC} provides a description of short curves.
In this paper, we prove the following ``convexity" property for the length 
of a curve along a \Teich geodesic. Let $g\co \R \to \T(S)$ be a geodesic 
in the \Teich space of $S$. For a curve $\alpha$ on $S$, denote the hyperbolic 
length of the geodesic representative of $\alpha$ at $g(t)$ by $l_t$.

\begin{theorem} \label{thm:convex}
Assume $\alpha$ is balanced at $t_\alpha$ and $s \geq t_\alpha$ 
(respectively, $s \leq t_\alpha$). Then, for any $t \geq s$ 
($t \leq s$), we have
$$\frac 1{l_s} \succ \frac1{l_t} .$$
\end{theorem}

We also give the following estimate for the twisting parameter along a \Teich 
geodesic. Let $\np$ be the stable foliation of the geodesic $g$. The twisting 
parameter around a curve $\alpha$ at $g(t)$ is (roughly) the number 
of times that $\np$ twists around $\alpha$ relative to a curve
perpendicular to $\alpha$ in the hyperbolic metric of $g(t)$, and is
denoted by $tw_t^+$.
 
\begin{theorem} \label{thm:tw-plus}
There exists a constant $d_\alpha>0$ such that
$$
tw^+_t(\alpha) =
\frac {d_\alpha\,e^{-2\,(t-t_\alpha)}}{e^{2\,(t-t_\alpha)}+e^{-2\,(t-t_\alpha)}} 
\pm O(1/l_t) 
$$
\end{theorem}

\subsection*{Some notation} 
To simplify our presentation, we avoid keeping track of constants
that depend on the topology of the surface only. Instead, we use
the following notation: When two functions $f$ and $g$ are equal
up to additive constants, that is, when there exists a $C$ depending on the 
topology of $S$, such that  $g(x) -C \leq f(x) \leq g(x)+C$, we write 
$f(x) \eadd g(x)$. Similarly, $f(x) \gadd g(x)$ and $f(x) \ladd g(x)$ mean
that the inequalities are true up to an additive constant. When an inequality 
is true up to a multiplicative constant, we use symbols $\emul$, $\gmul$ 
and $\lmul$; and, when it is true up to an additive constant and a 
multiplicative constant, 
we use symbols $\asymp$, $\prec$ and $\succ$. For example,
$f(x) \asymp g(x)$ means that there are constants $c$ and $C$, depending
on the topology of the surface only, such that
$$\frac 1c \, g(x) -C \leq f(x) \leq c\, g(x) + C.$$

\subsection*{Acknowledgments}
I would like to thank Yair Minsky for suggesting the statement of 
\thmref{thm:main} and Jason Behrstock for persuading me to
prove a more general version and suggesting the statement of 
\thmref{thm:general}. I would also like to thank Yair Minsky and
Young-Eun Choi for many helpful conversations.
 
\section{Preliminaries} \label{sec:prel}

\subsection{Curves and markings}
By a {\it curve} in $S$ we mean a non-trivial, non-peripheral, simple closed 
curve in $S$. The free homotopy class of a curve $\alpha$ is denoted by
$[\alpha]$. By an {\it essential arc} $\omega$ we mean a simple arc, 
with endpoints on the boundary of $S$, that cannot be pushed to the 
boundary of $S$. In case $S$ is not an annulus, $[\omega]$ represents 
the homotopy class of $\omega$ relative to the boundary
of $S$. When $S$ is an annulus, $[\omega]$ is defined to be  the homotopy 
class of $\omega$ relative to the endpoints of $\omega$. 

Define $\CC(S)$ to be the set of 
all homotopy classes of curves and essential arcs on the surface $S$. 
To simplify notation, we often write $\alpha \in \CC(S)$ instead of
$[\alpha] \in \CC(S)$. Define a distance on $\CC(S)$ as follows: 
For $\alpha, \beta \in \CC(S)$, define
$d_S(\alpha,\beta)$ to be equal to one if $\alpha \not = \beta$ 
and if $\alpha$ and $\beta$ can be represented by disjoint curves
or arcs. Let the metric on $\CC(S)$ be the maximal metric having the 
above property, i.e., $d_S(\alpha, \beta)= n$ if 
$\alpha=\gamma_0, \gamma_1, \ldots, \gamma_n=\beta$ is the 
shortest sequence of curves or arcs on $S$ such that, for
$i=1, \ldots, n$, $\gamma_{i-1}$ is distance one from $\gamma_i$. 
(See \cite{minsky:CCI}.) 

Let $\{ \alpha_1, \ldots , \alpha_m \}$ be a pants 
decomposition of $S$. A {\it marking} on $S$ is a set 
$\mu=\{ (\alpha_1, \beta_1), \ldots , (\alpha_m, \beta_m) \}$ such that the 
curve $\beta_i$ is disjoint from $\alpha_j$, for $i \not = j$, 
and intersects $\alpha_i$ once (twice) if the surface filled by $\alpha_i$ 
and $\beta_i$ is a once-punctured torus (four-times-punctured sphere).
The $\alpha_i$ are called the base curves of $\mu$.  For every 
$i$, $\beta_i$ is called the transverse curve to $\alpha_i$ in
$\mu$.
When the distinction between the base curves and the transverse curves
is not important, we represent a marking as a set of curves
$\{ \beta_1, \ldots, \beta_n\}$ including all the base curves and the 
transverse curves. Denote the space of all markings on $S$ by $\M(S)$ 
(see \cite{minsky:CCII}.)

\subsection{Subsurface intersection and subsurface distance}
Let $\nu$ be a subset of $\CC(S)$ (e.g., curves appearing in a marking) 
or a singular foliation on $S$, and let $Y$ be a subsurface of $S$. We define 
the {\it projection} of $\nu$ to the subsurface $Y$ as follows: Let
$$f:\bar S \to S$$
be a regular covering of $S$ such that $f^*(\pi_1(\bar S))$ is conjugate to
$\pi_1(Y)$ (the $Y$--cover of $S$). Since $S$ admits a hyperbolic 
metric, $\bar S$ has a well-defined boundary at infinity. Let $\bar \nu$ be 
the lift of $\nu$ to $\bar S$. Components of $\bar \nu$ that are 
essential arcs or curves on $\bar S$, if any, form a subset of $\CC(\bar S)$.
The surface $\bar S$ is homeomorphic to $Y$. We call the corresponding 
subset of $\CC(Y)$ the projection of $\nu$ to $Y$ and will denote it by $\nu_Y$. 
If there are no essential arcs or curves in $\bar \nu$, $\nu_Y$ 
is the empty set; otherwise we say that {\it $\nu$ intersects $Y$ essentially}.
This projection depends on the homotopy class of elements of $\nu$ only.

Let $\nu$ and  $\nu'$ be subsets of $\CC(S)$ or singular foliations on $S$
that intersect a subsurface $Y$ essentially. We define the
{\it $Y$--intersection} ({\it $Y$--distance}\/) between $\nu$ and $\nu'$ 
to be the maximum geometric intersection number in $Y$ 
(maximum distance in $\CC(Y)$) between the elements of projections 
$\nu_Y$ and $\nu_Y'$ and denote it by
$$
i_Y(\nu,\nu') \qquad \big(\text{respectively,} \quad d_Y(\nu,\nu') \big).
$$
If $Y$ is an annulus whose core is the curve $\alpha$, then we also denote $i_Y(\nu,\nu')$ and $d_Y(\nu,\nu')$ by $i_\alpha(\nu,\nu')$ and $d_\alpha(\nu,\nu')$, respectively. The following lemma is well known.

\begin{lemma} \label{lem:dist-int}
Let $Y$, $\nu$ and $\nu'$ be as above. 
\begin{enumerate}
\item If $Y$ is not an annulus, then
$$d_Y(\nu,\nu') \prec \log i_Y(\nu,\nu').$$
\item For a curve $\alpha$,
$$d_\alpha(\nu,\nu') \asymp i_\alpha(\nu,\nu').$$
\end{enumerate}
\end{lemma}

\subsection{Quadratic differentials} \label{sec:quad}
Let $q$ be a meromorphic quadratic differential of area one on $S$.
(See \cite{gardiner:QT} for definition and details.)
We assume that $q$ has a discrete set of finite critical points 
(i.e., critical points of $q$ are either zeroes or poles of order 1). 
Corresponding to $q$, there are two singular measured foliations called 
the horizontal and the vertical foliations, which we denote by $\np$ and 
$\nm$. We call the singular Euclidean metric
$|q|$ the {\it $q$--metric} on $S$. For a curve $\alpha$ in $S$, the 
$q$--geodesic representative of $\alpha$ exists and is unique except for 
the case where it is one of the continuous family of closed geodesics 
in a flat annulus, which we refer to as the flat annulus corresponding to 
$\alpha$. (Some difficulties aries when $q$ has poles of order 1.
See \cite{me:SC} for precise definitions and discussion.)
We denote the $q$--length of $\alpha$ by $l_q(\alpha)$, 
the horizontal length of $\alpha$ by $h_q(\alpha)$ and the vertical length 
of $\alpha$ by $v_q(\alpha)$. We also denote the $q$--length, the horizontal
length and the vertical length of the $q$--geodesic representative of 
$\alpha$, by $l_q([\alpha])$, $h_q([\alpha])$ and
$v_q([\alpha])$, respectively. In general, for any metric $\tau$, 
$l_\tau(\alpha)$ represents the $\tau$--length of $\alpha$ and 
$l_\tau([\alpha])$ represents the $\tau$--length of the $\tau$-geodesic
representative of $\alpha$.

\subsection{Regular and primitive annuli in $q$} 
Let $Y$ be a subsurface of $S$ and $\gamma$ be a boundary 
component of $Y$.\footnote{We always assume that curves are
piecewise smooth.}  The curvature of $\gamma$ with respect to $Y$, 
$\kappa_Y(\gamma)$, is well defined as a measure with atoms at the 
corners. We choose the sign to be positive when the acceleration vector 
points into $Y$. If $\gamma$ is curved non-negatively (or
non-positively) with 
respect to $Y$ at every point, we say it is {\it monotonically curved} with 
respect to $Y$. Let $A$ be an open annulus in $S$ with boundaries 
$\gamma_0$ and $\gamma_1$. Suppose both boundaries are monotonically
curved with respect to $A$ and $\kappa_A(\gamma_0) \leq 0$.
Further, suppose that the boundaries are equidistant from each other, and 
the interior of $A$ contains no zeroes. We call $A$ a {\it primitive} annulus 
and write $\kappa(A)=-\kappa_A(\gamma_0)$. If $\kappa(A) > 0$,
we call $A$ {\it expanding}  and say that $\gamma_0$ is the inner 
boundary and $\gamma_1$ is the outer boundary. When $\kappa(A)=0$, 
$A$ is a {\it flat} annulus and is foliated by closed Euclidean geodesics 
homotopic to the boundaries. The following lemma is useful for computing
the modulus of a primitive annulus. 

\begin{lemma}[{\cite[Lemma 3.6]{me:SC}}] \label{lem:p-dist}
Let $A$ and $\gamma_0$ be as above, and let $d$ be the distance 
between the boundaries of $A$. Then
\begin{equation*}
\left\{ \begin{array}{ll} \displaystyle
\kappa \, \Mod(A) \asymp \log \left( \frac{d}{ l_q(\gamma_0)} \right) &  
                                        \text{if} \: \: \kappa(A) > 0 \\ & \\
\Mod(A) \, l_q(\gamma_0) = d  & \text{if} \: \: \kappa(A) = 0
\end{array} \right. .
\end{equation*}
\end{lemma}

Minsky has shown that every annulus of large modulus contains
a primitive annulus with comparable modulus.

\begin{theorem}[Minsky {\cite[Theorem 4.6]{minsky:HM}}] \label{thm:primitive}
There exists an $\ep_0>0$ such that, for a  curve $\alpha$ in $S$, if
$l_\sigma([\alpha]) \leq \ep_0$, then there exists a primitive annulus $A$ 
such that 
$$\frac 1{l_\sigma([\alpha])} \asymp \Mod(A).$$
\end{theorem}

Throughout this paper, $\ep_0$ is a fixed constant smaller than
the Margulis constant, such that the above theorem and \thmref{thm:product}
are true.

\subsection{Product regions in \Teich space}
The \Teich space of $S$, $\T(S)$, is the space of conformal structures
on $S$ up to isotopy. The Teichm\"uller distance between two points
$\sigma_1$ and $\sigma_2$ is defined as
$$d_\T(\sigma_1,\sigma_2) = \frac 12 K(\sigma_1,\sigma_2),$$
where $K(\sigma_1,\sigma_2)$ is the smallest quasi-conformal dilatation
of a homeomorphism from $\sigma_1$ to $\sigma_2$. 
Let $\Gamma$ be a system of disjoint curves on $S$, and let 
$\thin_\ep(\Gamma)$ denote the set of all $\sigma \in \T(S)$ such that,
 for all $\gamma \in \Gamma$, the length of $\gamma$ in $\sigma$,
$l_\sigma(\gamma)$, is less than or equal to $\ep$. Let 
$\T_\Gamma$ denote the product space
$$\T(S \setminus \Gamma) \times 
\prod_{\gamma \in \Gamma} \hyp_\gamma,$$
where $S \setminus \Gamma$ is considered as a punctured space
and each $\hyp_\gamma$ is a copy of the hyperbolic plane. Endow 
$\T_\Gamma$ with the sup metric. Minsky has shown, for small
enough $\ep$, that $\thin_\ep(\Gamma)$ has a product structure.

\begin{theorem}[Minsky \cite{minsky:PR}] \label{thm:product}
The Fenchel-Nielsen coordinates on $\T(S)$ give rise to
a natural homeomorphism $\pi \co \T(S) \to \T_\Gamma$. 
There exists an $\ep_0> 0$ sufficiently small that this homeomorphism
restricted to $\thin_{\ep_0}(\Gamma)$ distorts distances by a bounded
additive amount.
\end{theorem}

Note that $\T(S \setminus \Gamma) = \prod_Y \T(Y)$, where the product
is over all connected components $Y$ of $S \setminus \Gamma$.
Let $\pi_0$ denote the component of $\pi$ mapping to 
$\T(S \setminus \Gamma)$, let $\pi_Y$ denote the component mapping 
to $\T(Y)$, and, for $\gamma \in \Gamma$, let $\pi_\gamma$ denote the 
component mapping to $\hyp_\gamma$. 
For the rest of the paper, we fix $L_0 >0$ \label{page:L} 
such that, for a hyperbolic metric $\sigma$ on $S$, if 
$l_\sigma(\alpha) \geq \ep_0$, then there exists a curve 
$\beta$ intersecting $\alpha$ with $l_\sigma(\beta) \leq L_0$. 

\section{Behavior of a Geodesic in the Thin Part of \Teich Space}
\noindent
In this section, we prove \thmref{thm:convex}, restated as
\thmref{thm:local}, and study how the combinatorics of short
markings changes along a \Teich geodesic. We show that, for every curve 
$\alpha$ in $S$, there exists a connected interval where $\alpha$ is ``short" 
(\corref{cor:interval}), and the projections of the short markings to a subsurface
can only change while all the boundaries of that subsurface are short
(\propref{prop:progress}). This is an essential component of the proof of
the main theorem.

\subsection{\Teich geodesics}
For $t \in \R$, let $q_t$ be the quadratic differential obtained from $q$ by 
scaling its horizontal foliation by a factor of $e^t$, and its vertical foliation 
by a factor of $e^{-t}$. Define $g(t)$ to be the conformal structure 
corresponding to $q_t$. Then $g \co \R \to \T(S)$ is a geodesic in 
$\T(S)$ parametrized by arc length. For a curve $\alpha$ in $S$, the horizontal 
and vertical lengths of $\alpha$ vary with time as follows: 
\begin{equation} \label{eq:hv-length}
h_{q_t}(\alpha) = h_q(\alpha) \, e^{-t} \qquad {\rm and} \qquad
v_{q_t}(\alpha) = v_q(\alpha) \, e^t.
\end{equation}
We say $\alpha$ is {\it balanced}, {\it mostly horizontal} or {\it mostly vertical} at time 
$t$ if, respectively, $v_t([\alpha])=h_t([\alpha])$, $v_t([\alpha]) \leq h_t([\alpha])$ 
or $v_t([\alpha]) \geq h_t([\alpha])$.

\subsection{Hyperbolic length along a geodesic}
The behavior of the hyperbolic length of a curve along
a \Teich geodesic is somewhat mysterious. For the Weil-Petersson 
metric on $\T(S)$, the hyperbolic length of a curve along a geodesic
is a convex function of time. In the \Teich metric, the quadratic
differential length of a curve is also convex. The following result is a weaker
but analogous statement. It roughly states that a curve assumes its shortest 
length when it is balanced and the length is ``non-decreasing" as one 
moves away in either direction.  Let $\sigma_t$ denote the hyperbolic metric on 
$g(t)$.

\begin{theorem} \label{thm:local}
Let $g$ be a geodesic in $\T(S)$ and $\alpha$ be a curve in 
$S$. Assume $\alpha$ is balanced at $t_\alpha$ and $s \geq t_\alpha$ 
(respectively, $s \leq t_\alpha$). Then, for any $t \geq s$ 
($t \leq s$), we have
\begin{equation} \label{eq:local}
\frac 1{l_{\sigma_s}([\alpha])} \succ \frac 1{l_{\sigma_t}([\alpha])}.
\end{equation}
\end{theorem}

\begin{remark}
The above inequality has no content if both $l_{\sigma_s}([\alpha])$ 
and $l_{\sigma_t}([\alpha])$ are large, because both quantities are within 
the additive error. However, if $l_{\sigma_s}([\alpha])$ is large, then
\eqref{eq:local} implies that $l_{\sigma_t}([\alpha])$ is bounded below
for all $t \geq s$.
\end{remark}

\begin{proof}
Let $F_t$ be the flat 
annulus corresponding to $\alpha$ in $q_t$. The modulus of $F_t$ is
maximum at $t_\alpha$, and, for $t \in \R$, 
\begin{equation} \label{eq:F}
\Mod(F_t) \asymp \Mod(F_{t_0}) \, e^{-2 \, |t-t_\alpha|}. 
\end{equation}
Let $A_t$ be as in \thmref{thm:primitive} for hyperbolic metric $\sigma_t$,
quadratic differential $q_t$, and curve $\alpha$ (if $l_t(\alpha) \geq \ep_0$,
there is nothing to prove). If $A_t$ is flat, then
\begin{align*}
\frac 1{l_{\sigma_s}([\alpha])}  
    &\asymp \frac 1{\Ext_{\sigma_s}([\alpha])} \tag{\cite{maskit:HE}}\\
    & \succ  \Mod(F_s)        \tag{by definition of $\Ext_{\sigma_s}(\alpha)$}\\
    & \succ  \Mod(F_t)         \tag{Equation \eqref{eq:F}}\\
    &\geq \Mod(A_t)            \tag{$A_t \subset F_t$}\\
    &\asymp \frac 1{l_{\sigma_t}([\alpha])}.   \tag{\thmref{thm:primitive}}
\end{align*}

Assume $A_t$ is not flat. Let $d$ be the distance between the 
boundary components of $A_t$ and $l$ be the length of
the inner boundary of $A_t$. Let $\beta$ be a curve intersecting
$\alpha$ whose hyperbolic length at $s$ is less than $L$, for some $L$ 
such that $e^L \asymp \frac 1{l_{\sigma_s}([\alpha])}.$
Using the ``collar lemma" (\thmref{thm:q-collar}), we have
\begin{equation} \label{eq:ba-ratio}
\frac 1{l_{\sigma_s}([\alpha])} \succ 
\log \frac {l_{q_s}([\beta])}{l_{q_s}([\alpha])}.
\end{equation}
But $\alpha$ is mostly vertical at $s$; therefore, for $t \geq s$,
$$l_{q_t}([\alpha]) \asymp l_{q_{t_\alpha}}([\alpha])\, e^{t-s}.$$
The quadratic differential length of any curve grows at most
exponentially; that is, for $t \geq s$,
$$l_{q_t}([\beta]) \prec l_{q_{t_\alpha}}([\beta])\, e^{t-s}.$$
Therefore,
\begin{equation}\label{eq:ratio-inc}
\frac {l_{q_s}([\beta])}{l_{q_s}([\alpha])} \geq 
\frac {l_{q_t}([\beta])}{l_{q_t}([\alpha])}.
\end{equation}
We also have $l_{q_t}([\beta]) \geq d$ ($\beta$ has to cross $A_t$)
and $l_{q_t}([\alpha]) \leq l$ ($\alpha$ and the inner boundary of $A_t$
are homotopic). Therefore, 

\begin{align*}
\frac 1{l_{\sigma_s}([\alpha])} 
   & \succ \log  \frac {l_{q_s}([\beta])}{l_{q_s}([\alpha])}
                                                              \tag{Equation \eqref{eq:ba-ratio}} \\
  & \geq  \log \frac {l_{q_t}([\beta])}{l_{q_t}([\alpha])} 
                                                              \tag{Equation \eqref{eq:ratio-inc}} \\
  & \geq  \log \frac dl                                \\
  & \succ \Mod(A_t)                                \tag{\lemref{lem:p-dist}}\\
  & \asymp \frac 1{l_{\sigma_t}([\alpha])}. \tag{\thmref{thm:primitive}}
\end{align*}
\end{proof}

\begin{corollary} \label{cor:interval}
There exists  $\ep_1$ such that, for any geodesic in the \Teich space and 
any curve $\alpha$ in $S$, there exists a connected (perhaps empty) 
interval $I_\alpha$ such that
\begin{enumerate}
\item for $t \in I_\alpha$, $l_{\sigma_t}([\alpha]) \leq \ep_0$, and
\item for $t \not \in I_\alpha$, $l_{\sigma_t}([\alpha]) \geq \ep_1$.
\end{enumerate}
\end{corollary}
The intersection of connected intervals is a connected interval
(or an empty set). Therefore, a similar statement is also true
for subsurfaces.

\begin{corollary} \label{cor:Y-interval}
Let $\ep_0$, $\ep_1$ and $g$ be as above. For every
subsurface $Y$, there exists a 
connected interval $I_Y $ such that
\begin{enumerate}
\item for $t \in I_Y$, the hyperbolic lengths of all boundary components of
$Y$ at $\sigma_t$ are less than or equal to $\ep_0$, and 
\item for $t \not \in I_Y$, there exists a boundary component of $Y$
whose hyperbolic length at $\sigma_t$ is greater than
or equal to $\ep_1$.
\end{enumerate}
\end{corollary}

\subsection{A lower bound for distance in the \Teich space}
Our main theorem describes how the distance between two points
in \Teich space can be estimated by measuring the combinatorial complexity
of curves of bounded size. Here we show that, if two curves of bounded
length in $\sigma_1$ and $\sigma_2$ intersect each other a large
number of times, then $\sigma_1$ and $\sigma_2$ are far apart
in $\T(S)$. 

First we recall some properties of the extremal length.
Let $\Ext_\sigma(\alpha)$ denote the extremal length of $\alpha$ in $\sigma$. 
Minsky has shown (see \cite{minsky:TG}) that, for curves $\alpha$ and $\beta$ 
in $S$, and $\sigma \in \T(S)$,
\begin{equation} \label{eq:int-bound}
\Ext_\sigma(\alpha) \, \Ext_\sigma(\beta) \geq i_S(\alpha, \beta)^2.
\end{equation}
Kerckhoff's theorem (see \cite{kerckhoff:AG}) states that, for points
$\sigma_1$ and $\sigma_2$ in $\T(S)$,
\begin{equation} \label{eq:ext-ratio}
K(\sigma_1,\sigma_2) = 
\sup_\alpha \frac {\Ext_{\sigma_1}(\alpha)}{\Ext_{\sigma_2}(\alpha)},
\end{equation}
where the sup is over all curves on $S$. 
We also know (see \cite{maskit:HE}) that, if the hyperbolic length of 
$\alpha$ is short (say, $l_\sigma(\alpha) \leq L_0$), then 
\begin{equation} \label{eq:ext-hyp}
l_\sigma(\alpha) \asymp \Ext_\sigma(\alpha).
\end{equation} 

\begin{proposition} \label{prop:int-dist}
Assume, for some $\sigma_1, \sigma_2 \in \T(S)$ and curves $\alpha$ and 
$\beta$ in $S$, that $l_{\sigma_1}(\alpha)\leq L_0$ and 
$l_{\sigma_2}(\beta) \leq L_0$. Then
$$d_\T(\sigma_1,\sigma_2) \succ \log i_S(\alpha,\beta).$$
\end{proposition}
\begin{proof}
We have:
\begin{align*}
i_S(\alpha,\beta)^2 
 &\leq \Ext_{\sigma_1}(\alpha) \, \Ext_{\sigma_1}(\beta) 
     \tag{Equation \eqref{eq:int-bound}}\\
 &\leq \Ext_{\sigma_1}(\alpha) \, \Ext_{\sigma_2}(\beta) \, 
    K(\sigma_1,\sigma_2) \tag{Equation \eqref{eq:ext-ratio}}\\
 & \asymp L_0^2 \; K(\sigma_1,\sigma_2) \tag{Equation \eqref{eq:ext-hyp}}
\end{align*}
Note that $L_0$ is a fixed constant depending on $S$ only. 
By taking the logarithm of both sides, we obtain the desired
inequality.
\end{proof}

\subsection{Combinatorics of short markings along a \Teich geodesic.}
For $t \in \R$, let $\mu_t$ be the {\it shortest marking} \label{page:short}
in $\sigma_t$, constructed as follows. Let $\alpha_1$ be the shortest curve 
in $S$ and $\alpha_2$ be the shortest curve disjoint from $\alpha_1$, and 
so on, to form a pants decomposition of $S$. Then, let the transverse
curve $\beta_i$ be the shortest curve intersecting $\alpha_i$ and
disjoint from $\alpha_j$, $i \not = j$.\footnote{There may be 
finitely many such markings.}
\propref{prop:progress} states that the projection of these markings
to a subsurface $Y$ stays in a bounded neighborhood in $\CC(Y)$ while 
the geodesic is outside of the thin part of $\T(S)$ corresponding to $Y$. 
The proof makes an essential use of the following theorem.

\begin{theorem}[{\cite[Theorem 5.5]{me:SC}}] \label{thm:if}
Let $\alpha$ be a curve in $S$, $\beta$ be a transverse curve to
$\alpha$, and $Y$ be a component of $S \setminus \alpha$
($Y$ is allowed to be an annulus). 
Assume $l_{\sigma_t}([\beta]) \leq L$. We have:
\begin{enumerate}
\item If $\alpha$ is mostly vertical, then
$$  i_Y(\beta, \np) \prec D_L.$$
\item If $\alpha$ is mostly horizontal, then
$$  i_Y(\beta, \nm) \prec D_L.$$
\end{enumerate}
Here, $D_L$ is a constant depending on $L$, with $\log D_L\asymp e^L$.
\end{theorem}
 
\begin{proposition} \label{prop:progress}
If $[r,s] \cap I_Y = \emptyset$, then 
$$d_Y(\mu_r,\mu_s) = O(1).$$
\end{proposition}
\begin{proof}
Let $L_1$ be such that every curve of length larger than $\ep_1$ in a 
hyperbolic surface with geodesic boundary has a transverse curve of 
length less than $L_1$.  For $t \in [r,s]$, there exists a boundary
component $\gamma_t$ of $Y$ whose $\sigma_t$--length is larger than 
$\ep_1$. Therefore, the marking $\mu_t$ contains a curve $\alpha_t$ 
with $l_{\sigma_t}(\alpha_t) \leq L_1$ that intersects $Y$ nontrivially. 
The projection of $\mu_t$  to $Y$ has bounded diameter. Therefore it is 
sufficient to prove $d_Y(\alpha_r,\alpha_s)=O(1)$.

The curve $\gamma_t$ is either mostly horizontal or mostly vertical at time $t$. 
The set of times at which $Y$ has a boundary component of length 
larger than or equal to $\ep_1$ which is mostly horizontal
(or mostly vertical) is closed. Therefore, either
\begin{enumerate}
\item $\gamma_r$ and $\gamma_s$ are both mostly horizontal or
both mostly vertical, or
\item for some $t \in [r,s]$, there are two curves $\gamma_t$ and 
$\gamma_t'$ whose lengths at $\sigma_t$ are larger than or equal to 
$\ep_1$, and one is mostly horizontal and the other is mostly vertical 
(possibly $\gamma_t = \gamma_t'$ and $\gamma_t$ is balanced).
\end{enumerate}

\noindent Case 1: 
If $\gamma_r$ and $\gamma_s$ are mostly vertical, \thmref{thm:if} implies that 
$$i_Y(\alpha_r, \np) \prec D_{L_1} \quad\text{and}\quad
i_Y(\alpha_s, \np) \prec D_{L_1} .$$
Therefore, using \lemref{lem:dist-int},
$$d_Y(\alpha_s,\np) \prec \log i_Y(\alpha_s, \np) \prec \log D_{L_1}.$$
Similarly, $d_Y(\alpha_r,\np) = O(1)$. This implies that
$d_Y(\alpha_r,\alpha_s)=O(1)$. The proof is similar if 
$\gamma_r$ and $\gamma_s$ are both mostly horizontal. \smallskip

\noindent Case 2: Assume (without loss of generality) that $\gamma_t$ is
mostly horizontal and $\gamma_t'$ is mostly vertical. Let $\alpha_t$ and 
$\alpha_t'$ be the corresponding transverse curves in $\mu_t$ of length 
less than $L_1$. By the above argument,
$$d_Y(\alpha_t,\nm) =O(1) \qquad \text{and} \qquad
d_Y(\alpha_t',\np) = O(1).$$
But the extremal lengths of $\alpha_t$ and $\alpha_t'$ are bounded 
by a constant depending on $L_1$. Equation \eqref{eq:int-bound} 
implies that $i_S(\alpha_t,\alpha_t') = O(1)$, and,
by  \lemref{lem:dist-int}, $d_Y(\alpha_t,\alpha_t') = O(1)$. Therefore,
\begin{equation} \label{eq:npnm}
d_Y(\np,\nm) =O(1).
\end{equation}
Again, as above, the projection of each of $\alpha_s$ and $\alpha_r$
to $Y$ is close to the projection of either $\np$ or $\nm$ to $Y$.
Thus, \eqref{eq:npnm} and the triangle inequality for $d_Y$ imply that 
\begin{equation*}
d_Y(\alpha_r,\alpha_s)=O(1).  \qedhere 
\end{equation*}
\end{proof}

\begin{corollary} \label{cor:NC-in-thick}
If $I_Y=[c,d] \subset [a,b]$, then 
$$d_Y(\mu_a,\mu_b) \asymp d_Y(\mu_c,\mu_d).$$
\end{corollary}

\section{Twisting in the Hyperbolic Metric vs. Twisting in the 
Quadratic Differential Metric} \label{sec:twisting}
\noindent
Let $\alpha$ be a curve in $S$. Having a metric in $S$ enables us
to define a twisting parameter for curves that cross $\alpha$.
This, roughly speaking, is the number of times that a given curve
twists around $\alpha$ in comparison with an arc that is perpendicular
to the geodesic representative of $\alpha$. In this section we define
a twisting parameter for $\np$ and $\nm$ using metrics given by
$q$ and $\sigma$, and we study how these two quantities are
related. We use this to prove \thmref{thm:tw-plus} at the end of this seection.

Let $\bar S$ be the annular cover of $S$ with respect to $\alpha$.
Let $\bar q$, $\npb$ and $\nmb$ be the 
lifts of $q$, $\np$ and $\nm$ to $\bar S$, respectively, and $\bar \beta_q$ 
be a geodesic arc connecting the boundaries of $\bar S$ that is perpendicular 
(in $\bar q$) to the geodesic representative of the core of $\bar S$, $\bar \alpha$.  
We define the twisting parameter of $\np$ around $\alpha$ in $q$
to be the maximum intersection number of a leaf of $\npb$ and 
$\bar \beta_q$, and we denote it by $tw_q(\np,\alpha)$. When it is
clear what $\alpha$ is, we denote this by $tw_q^+$.  The twisting 
parameter $tw_q^-$ of $\nm$ around $\alpha$ in $q$ is defined
similarly. Note that the maximum intersection number is at least
one, that is, $tw_q^\pm$ are positive integers. 

\begin{figure}[htbp]
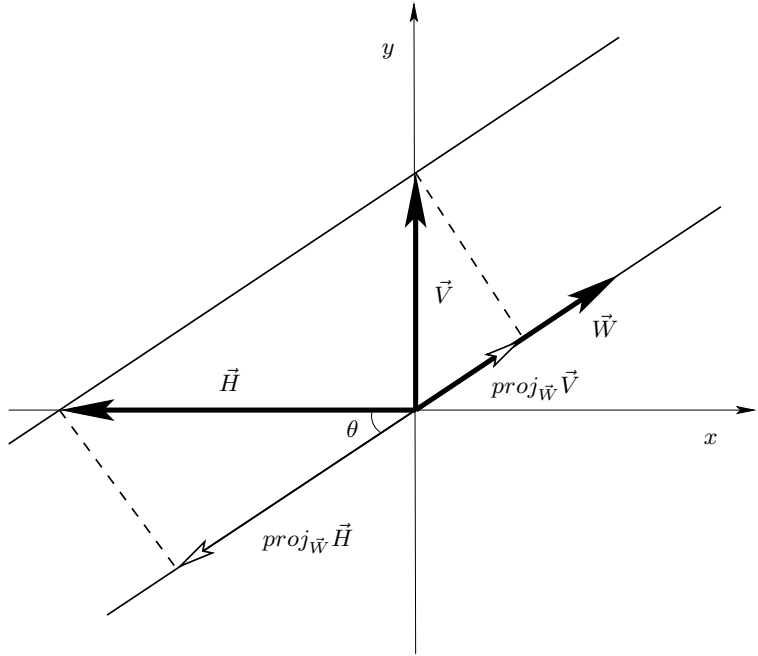

\begin{center}
\input annulus.pstex_t
\caption{\bf The universal cover of $F$}
\label{fig:annulus}
\end{center}
\end{figure}

Let $F$ be the flat annulus in $q$ corresponding to $\alpha$ and
let $\beta_q$ be an arc connecting the boundaries of $F$ that
is perpendicular to the boundaries of $F$. The intersection number
of the lift of a leaf of $\np$ with $\bar \beta_q$ is (up to small additive error) 
equal to the intersection number of the restriction of this leaf to $F$
with $\beta_q$.  Therefore, to compute $tw_q^\pm$, it is sufficient
to understand the picture in $F$. Consider an isometric embedding of the universal cover of $F$ in $\R^2$ 
such that the leaves of horizontal foliations are parallel to the $x$--axis 
and the leaves of vertical foliations are parallel to the $y$-axis 
(see \figref{fig:annulus}). 

Let $\vec W$ be the vector representing the translation that
generates the deck translation group.  Let $\vec H$ be the lift of a leaf of 
$\np$ passing through the origin 
and $\vec V$ be the same for $\nm$. From the above discussion, we have:
$$tw_q(\np,\alpha) \eadd \frac{\| \text{Proj}_{\vec W} \vec H \|}{\| \vec W \|}
\quad\text{and}\quad 
tw_q(\nm,\alpha) \eadd \frac{\| \text{Proj}_{\vec W} \vec V \|}{\| \vec W \|}. $$

Let $\theta$ be the angle between $\vec W$ and the $x$--axis.
It is easy to see, using similar triangles, that 
$$
\frac{\| \text{Proj}_{\vec W} \vec H \|}{\| \text{Proj}_{\vec W} \vec V \|} =
\frac {\sin^2 \theta}{\cos^2 \theta}.$$
We also have 
$\frac{h_q([\alpha])}{v_q([\alpha])} = \frac {\sin \theta}{\cos \theta}$. 
Therefore,
\begin{equation} \label{eq:q-twisting}
\frac{tw_q^+}{tw_q^-} \eadd \frac{h_q([\alpha])^2}{v_q([\alpha])^2}. 
\end{equation}

This is a very useful equation that allows us to compute the $q$--twisting
parameter of horizontal and vertical foliations around $\alpha$ along a \Teich 
geodesic (see equation \eqref{eq:qt-twist}).

We define the twisting parameter for a hyperbolic metric
as follows. Let $\beta_\sigma$ be the shortest transverse curve to 
$\alpha$ in the hyperbolic metric $\sigma$. Define
$$tw_\sigma^+ = i(\np, \beta_\sigma) \quad\text{and}\quad 
tw_\sigma^- = i(\nm, \beta_\sigma).$$
We would like to prove a statement similar to equation \eqref{eq:q-twisting}
for $\sigma$-twisting parameters.  However, giving good estimates
for $tw_\sigma^\pm$ is difficult when $\alpha$ is very short. The errors 
in our estimates get larger as $l_\sigma(\alpha)$ gets smaller.

Let $\bar \beta_\sigma$ be the lift of $\beta_\sigma$ to $\bar S$
whose end points are in different boundary components of $\bar S$.
Our strategy is to relate $q$-- and $\sigma$--twisting parameters 
by providing an upper bound for $i(\bar \beta_q, \bar \beta_\sigma)$.

\begin{lemma} \label{lem:n}
If $i (\bar \beta_q,\bar \beta_\sigma) =n$, then
$$\Ext_{\sigma}(\beta_\sigma) \gmul  n ^2 \, l_\sigma(\alpha).$$
\end{lemma}

\begin{proof}
By definition of the extremal length, for any metric $\tau$ on $S$ in the 
conformal class of $\sigma$, 
$$
\Ext_\sigma(\beta_\sigma) \geq \frac {l_\tau(\beta_\sigma)^2}{\area_\tau(S)}.
$$
To find a lower bound for $\Ext_\sigma(\beta_\sigma)$, we need to find an
appropriate metric $\tau$. First we establish some notation. 
Let $A$ be the largest regular neighborhood of $F$ that is still an 
annulus. Denote the boundary components of $A$ by $\alpha_0$ and 
$\alpha_c$, where $c$ is the $q$--distance between the boundaries of $A$. 
For $t \in (0,c)$, let $\alpha_t$ be a curve in $A$ that is equidistant from
a $q$--geodesic representative of $\alpha$ and whose 
$q$--distance from $\alpha_0$ is $t$. These curves give a foliation
of $A$ into curves in the homotopy class of $\alpha$. There is a 
subinterval $[a,b]$ of $[0,c]$ such that, for $t\in[a,b]$, $\alpha_t$
is a $q$--geodesic representative of $\alpha$. This gives a division of
$A$ into three pieces, the flat annulus $F$ containing all $\alpha_t$, 
$t\in[a,b]$, and two expanding annuli $A_1$ and $A_2$ on the sides.
\thmref{thm:primitive} implies that $\Mod(A) \asymp \frac 1{l_\sigma(\alpha)}$.
Using \lemref{lem:p-dist}, we have
\begin{equation} \label{eq:abc}
\frac 1{l_\sigma(\alpha)}\emul
\log \frac{a}{l_q([\alpha])}  + \frac{(b-a)}{l_q([\alpha])} 
+ \log \frac{(c-b)}{l_q([\alpha])}.
\end{equation}
 
As $t$ changes in the interval $[b,c]$, the length of $\alpha_t$
increases.  The rate of change is equal to the curvature
of $\alpha_t$, which is bounded above and below by constants
depending on the topology of $S$ only. A similar statement is true for
$A_1$ as well. Therefore,
\begin{equation} \label{eq:length}
l_q(\alpha_t) \emul \left\{ \begin{array}{ll}
l_q([\alpha]) + (a-t) & \text{if}\quad t\in[0,a]\\
l_q([\alpha]) & \text{if}\quad t\in[a,b]\\
l_q([\alpha]) + (t-b)  & \text{if}\quad t\in[b,c]
\end{array} \right. .
\end{equation}
Denote $l_q(\alpha_t)$ by $\lambda_t$.

Let $Z$ be the union of $A$; the $\lambda_0$--neighborhood, $N_0$, of
$\alpha_0$; and the $\lambda_c$--neighborhood, $N_c$, of $\alpha_c$. 
Define the metric $\tau$ in $S$ in the conformal class of $q$ as 
follows: if $x$ lies on a curve $\alpha_t$ in $A$, then we scale the $q$--metric
at $x$ by a factor of $\frac 1{\lambda_t}$; if $x$ is outside of $A$ and in
$N_0$, then  we scale the $q$--metric at $x$ by a factor of 
$\frac 1{\lambda_0}$; if $x$ is outside of $A$ and in $N_c$, then 
we scale the $q$--metric at $x$ by a factor of
$\frac 1{\lambda_c}$ (if $x$ is in both $N_0$ and 
$N_c$, then we scale the $q$--metric by a factor of 
$\max (\frac 1{\lambda_0}, \frac 1{\lambda_c})$); and, if $x$ is outside 
of $Z$, then we scale the $q$--metric at $x$ by a small enough factor so that the 
$\tau$--area of $S$ is comparable with the $\tau$--area of $Z$.
Note that $\area_q N_0 \lmul \lambda_0^2$
and $\area_q N_c \lmul \lambda_c^2$. We have
\begin{align*}
\area_\tau(S) & \emul \area_\tau(Z) \\
  & \leq \area_\tau N_0 + \area_\tau N_c + \area_\tau A\\
  & \emul 1+1+\int_0^c 1 \, . \, \frac {dt}{\lambda_t} \\
  & \emul 2 + \log \frac{a}{l_q([\alpha])}  + \frac{(b-a)}{l_q([\alpha])} 
     + \log \frac{(c-b)}{l_q([\alpha])} \tag{Equation \eqref{eq:length}} \\
  & \emul \frac 1{l_\sigma(\alpha)} \tag{Equation \eqref{eq:abc}}.
\end{align*}

Let  $\bar A$ be the lift of $A$ to $\bar S$ that is an annulus, 
and let $\bar \alpha_t$ be the lift of $\alpha_t$ that is in $\bar A$
(this is to ensure that $\bar \alpha_t$ is a closed curves not an infinite line).
Let $\bar \omega$ be a sub-arc of $\bar \beta_\sigma$ with end 
points in $\bar \beta_q$ that goes around $\bar S$ once, that is, if
$\bar \omega'$ is the sub-arc of $\bar \beta_q$ connecting the end points 
of $\bar \omega$, then $\bar \gamma= \bar \omega \cup \bar \omega'$ 
is a curve in the homotopy class of the core of $\bar S$. Let $\gamma$ 
be the projection of $\bar \gamma$ to $S$. Then $\gamma$ is in the 
homotopy class of $\alpha$
and therefore must intersect $A$ (otherwise, $A$ would not be maximal).
Hence, $\bar \gamma$ must intersect $\bar A$.  But $\bar \beta_q$
is perpendicular to $\bar \alpha_t$, and, once it exits $\bar A$, it never
returns. Therefore, $\bar \omega$ must intersect $\bar A$ as well. 

Let $\bar \alpha_s$ be an equidistant curve in $\bar A$ intersecting 
$\bar \omega$ that has the shortest $\bar q$--length . We claim that
$$l_{\bar q}(\bar \omega) \geq l_{\bar q}(\bar \alpha_t)= \lambda_t.$$
Assume $s > b$. The curve $\bar \alpha_s$ divides $\bar S$ into
two annuli. Let $B$ be the annulus that contains $\bar \alpha_c$.
For $t\in [b,s)$, the $\bar q$--length of $\bar \alpha_t$ is less
than the $\bar q$--length of $\alpha_s$. By assumption $\bar \alpha_s$
is the shortest equidistant curve intersecting $\bar \omega$,
therefore, $\bar \omega \subset B$. 

The curvature of $\bar \alpha_t$ with respect to $B$ is non-positive at
all points. Therefore, the closest-point projection from $B$ to 
$\bar \alpha_t$ is length-decreasing. But the end points of $\bar \omega$ project to 
the same point in $\bar \alpha_t$
(because $\bar \beta_q$ is perpendicular to $\bar \alpha_t$), and 
the projection covers $\bar \alpha_t$ completely. Therefore,
$l_{\bar q}(\bar \omega) \geq l_{\bar q}(\bar \alpha_t)$ in this case. 

A similar argument holds if $t < a$. If $t \in [a,b]$, then $\bar \omega$ 
could intersect $\bar \alpha_t$ transversally, but, in this case, $\bar \alpha_t$ is 
a $\bar q$--geodesic and the curvature of $\bar \alpha_t$ is non-positive with 
respect to both annuli in $\bar S \setminus \alpha_t$. Therefore, the claim 
is true in all cases.

Let $\omega$ be the projection of $\bar \omega$ to $S$. If $\omega$
exits $Z$, then its $\tau$--length is larger than the $\tau$--distance 
between $A$ and $\p Z$, which is equal to 1. Otherwise, 
$\omega \subset Z$. Then, at each point in $\omega$,
$\tau$ is obtained from $q$ by scaling by a factor of at least
$\frac 1{\lambda_t}$.  Therefore,
$$l_\tau(\omega) \geq \frac 1{\lambda_t} \, l_q(\omega) \geq 1.$$
There are $(n-1)$ arcs like $\bar \omega$, and they all
project down to different sub-arcs of $\beta_\sigma$. Therefore,
$$l_\tau(\beta_\sigma) \geq n.$$
This implies that
\begin{equation*}
\Ext_\sigma(\beta_\sigma) \geq \frac{l_\tau(\beta_\sigma)^2}{\area_\tau S} 
\gmul \frac{n^2}{1/l_\sigma(\alpha)} = n^2\,l_\sigma(\alpha). \qedhere
\end{equation*}
\end{proof}

\begin{corollary} \label{cor:twist-comparison}
For $\bar \beta_\sigma$ and $\bar \beta_q$ as before, we have
$$i(\bar \beta_\sigma,\bar \beta_q) \prec  \frac 1{l_\sigma(\alpha)}.$$
\end{corollary}

\begin{proof}
The curve $\beta_\sigma$ is the shortest (in $\sigma$) transverse 
curve to $\alpha$. Therefore, 
$\Ext_\sigma(\beta_\sigma) \asymp \frac 1{l_\sigma(\alpha)}$. 
Applying the previous theorem we get
$$  \frac 1{l_\sigma(\alpha)} \succ n^2 \, l_\sigma(\alpha),$$
which, using \lemref{lem:n}, implies the corollary.
\end{proof}

The following theorem is an immediate consequence of the definitions
of the twisting parameters and of \corref{cor:twist-comparison}.

\begin{theorem} \label{thm:twist-comparison}
The two twisting parameters are the same up to an additive error
comparable to $\frac 1{l_\sigma(\alpha)}$. That is,
$$tw_\sigma^\pm = tw_q^\pm \pm O(\frac 1{l_\sigma(\alpha)}).$$
\end{theorem}

\subsection{The twisting parameter along a \Teich geodesic.}
In this section, we give estimates for the twisting parameters of
$\npm$ around a curve $\alpha$ in $\sigma_t$.
Let $d=i_\alpha(\np,\nm)$. If $\alpha$ is not very short in $\sigma_t$,
say $l_{\sigma_t}(\alpha) \geq \ep_0$, 
then it has a transverse curve that is not longer than $L_0$.
\thmref{thm:if} implies that
\begin{equation} \label{eq:asymp}
\left\{ \begin{array}{ll}
\text{if $\alpha$ is mostly horizontal,}& 
tw_{\sigma_t}^+ \eadd d_\alpha \quad\text{and}\quad tw_{\sigma_t}^- \eadd 0\\
\text{if $\alpha$ is mostly vertical,}& 
tw_{\sigma_t}^+ \eadd 0 \quad\text{and}\quad tw_{\sigma_t}^- \eadd d_\alpha\\
\end{array} \right. .
\end{equation}
In general, we know that $tw_{q_t}^+ + tw_{q_t}^- \eadd d$. Assume
$\alpha$ is balanced at $t_\alpha$. Using Equations \eqref{eq:q-twisting} 
and \eqref{eq:hv-length}, we get
$$\frac {tw_{q_t}^+}{tw_{q_t}^-} = 
\frac {e^{ 2(t-t_\alpha)} \, h_{q_{t_\alpha}}(\alpha)^2}
        {e^{-2(t-t_\alpha)} \, v_{q_{t_\alpha}}(\alpha)^2}. $$
But $h_{q_{t_\alpha}}(\alpha) = v_{q_{t_\alpha}}(\alpha)$. Therefore,
\begin{equation} \label{eq:qt-twist}
tw_{q_t}^+ \eadd 
\frac{d_\alpha e^{2(t-t_\alpha)}}{e^{2(t-t_\alpha)}+e^{-2(t-t_\alpha)}} 
\quad\text{and}\quad
tw_{q_t}^- \eadd 
\frac{d_\alpha e^{-2(t-t-\alpha)}}{e^{2(t-t-\alpha)}+e^{-2(t-t_\alpha)}}.
\end{equation}
This and \thmref{thm:twist-comparison} prove \thmref{thm:tw-plus}.
The following theorem is a different statement for the same basic fact.

\begin{proposition} \label{prop:bounded-move}
Let $\sigma_t \in \T(S)$ and $\alpha$ be a curve in $S$ with
$l_{\sigma_t}(\alpha) \leq \ep_0$. Let $\sigma_t'$ be the point in $\T(S)$ 
obtained from $\sigma_t$ by twisting along $\alpha$ such that 
$$ tw_{\sigma_t'}^+ =  
    \frac{d-\alpha e^{-2(t-t_\alpha)}}{e^{2(t-t_\alpha)}+e^{-2(t-t_\alpha)}}.$$
Then $d_\T(\sigma_t,\sigma_t') = O(1)$.
\end{proposition}

\begin{proof}
Consider $\pi \co \T(S) \to \T(S \setminus \alpha) \times \hyp_\alpha$. 
We know that $\pi_0(\sigma_t)= \pi_0(\sigma_t')$ and
$$d_{\hyp_\alpha}(\pi_\alpha(\sigma_t), \pi_\alpha(\sigma_t')) \eadd
\log \big( l_{\sigma_t}(\alpha)\, (tw_{\sigma_t}^+ - tw_{\sigma_t'}^+) \big).$$
\thmref{thm:twist-comparison} implies that the $\sigma_t$--twisting
and the $q_t$--twisting parameters of $\np$ are equal up to an additive error
that is comparable with $\frac 1{l_{\sigma_t}(\alpha)}$. Therefore, the
right-hand side of the above equation is uniformly bounded. We have
\begin{equation*}
d_\T(\sigma_t, \sigma_t') \eadd 
d_{\hyp_\alpha}(\pi_\alpha(\sigma_t), \pi_\alpha(\sigma_t')) = O(1). \qedhere
\end{equation*}
\end{proof}

\section{Proof of the main theorem} \label{sec:proof}
\noindent
In this section we prove \thmref{thm:main}. In \secref{subsec:lower}, 
we show how a lower bound for the \Teich distance between two points
in $\T(S)$ can be obtain by the combinatorial complexity between
their short markings. In \secref{subsec;upper}, we give an upper bound for 
the distance between two points in the \Teich space by constructing a path 
in $\T(S)$ of length comparable with the estimate given in \thmref{thm:main}. 

\subsection{Lowers estimate} \label{subsec:lower} 
Let $g \co [a,b] \to \T(S)$ be the geodesic segment in the 
\Teich space connecting $\sigma_a$ to $\sigma_b$. Recall
that $\sigma_t$ is the hyperbolic metric of $g(t)$, and $\mu_t$
is the short-marking on $S$ corresponding to $\sigma_t$.

\begin{lemma} \label{lem:length}
Let $Y$ be a subsurface that is not an annulus and $I= I_Y \cap [a,b]$. Then
$$ |I| \succ d_Y(\mu_a,\mu_b).$$
\end{lemma}

\begin{proof}
Let $I=[c,d]$, $\tau_c=\pi_Y(\sigma_c)$ and $\tau_d= \pi_Y(\sigma_d)$ (see \thmref{thm:product}). Let $\eta_c$ and $\eta_d$ be the short-markings 
on $Y$ corresponding to $\tau_c$ and $\tau_d$, respectively. In fact, 
$\eta_c \subset \mu_c$ and $\eta_d \subset \mu_d$. We have
\begin{align*}
|I| &\succ d_{\T(Y)}(\tau_c,\tau_d), \tag{\thmref{thm:product}} \\
    &\succ \log i_Y(\eta_c,\eta_d) \tag{\propref{prop:int-dist}} \\
    &\succ  d_Y(\eta_c,\eta_d). \tag{\lemref{lem:dist-int}}
\end{align*}

But $d_Y(\eta_c,\eta_d) \asymp d_Y(\mu_c,\mu_d)$ (because
they have the same projections to $Y$). Also, by \propref{prop:progress},
we have 
$$d_Y(\mu_a,\mu_c) =O(1) \quad\text{and}\quad d_Y(\mu_d,\mu_b) =O(1).$$
This proves the lemma.
\end{proof}

A similar lemma is true when the subsurface is an annulus. The difference
is that, in \lemref{lem:length}, there is no restriction on the lengths
of the boundaries of $Y$; but, for the next lemma to be true, we have to assume
that $\alpha$ is not very short in $\sigma_a$ and $\sigma_b$. the proofs are
almost identical.

\begin{lemma} \label{lem:length-ann}
Let $\alpha$ be a curve in $S$ such that $l_{\sigma_a}(\alpha)\geq \ep_0$
and $l_{\sigma_b}(\alpha)\geq \ep_0$, and let $I=I_\alpha \cap [a,b]$. Then
$$ |I| \succ d_\alpha(\mu_a,\mu_b).$$
\end{lemma}

\begin{proof}
Since $\alpha$ is not short at either end, either $I_\alpha$ is disjoint 
from $[a,b]$ or it is a subset of $[a,b]$. If $I_\alpha \cap [a,b]=\emptyset$,
then \propref{prop:progress} implies the lemma. If
$I_\alpha = [c,d] \subset [a,b]$, then, by \corref{cor:NC-in-thick},
$$d_\alpha(\mu_a,\mu_b) \asymp d_\alpha(\mu_c,\mu_d).$$
Let $\beta_c$ and $\beta_d$ be curves transverse to $\alpha$
in markings $\mu_c$ and $\mu_d$, respectively. We have
$$i(\beta_c,\beta_d) = d_\alpha(\mu_c,\mu_d).$$
As in the previous lemma, using \thmref{thm:product} and  
\propref{prop:int-dist}, we have
$$|I_\alpha| \succ \log i(\beta_c,\beta_d).$$
The combination of the last three equations proves the lemma.
\end{proof}

The following proposition provides a lower bound for the \Teich distance 
between two points in the thick part of $\T(S)$.

\begin{proposition} \label{prop:lower}
Let $\sigma_1$, $\sigma_2$ be in the $\ep_0$--thick part of $\T(S)$ and $\mu_1$ and $\mu_2$ be the short-markings in $\sigma_1$ and 
$\sigma_2$, respectively. There exists a $k_0 > 0$ such that 
\begin{equation*}
d_\T(\sigma_1, \sigma_2 ) \succ \sum_Y \big[ d_Y(\mu_1, \mu_2) \big]_{k_0} 
+ \sum_\alpha \log \big[ d_\alpha(\mu_1,\mu_2) \big]_{k_0}.
\end{equation*}
\end{proposition}
\begin{proof}
Let $g\co [a,b] \to \T(S)$ be the geodesic segment connecting 
$\sigma_1$ and $\sigma_2$. Since the end points are in the thick
part of $\T(S)$, for every subsurface $Y$, $I_Y$ either is disjoint
from $[a,b]$ or is a subset of $[a,b]$. Let $k_0$ be a constant such that, if 
$d_Y(\mu,\eta) \geq k_0$, then $I_Y \subset [a,b]$ (see \propref{prop:int-dist}).
For $t \in I_Y$, the length of each boundary component of $Y$
is less than $\ep_0$. Therefore, there exists a constant $C$, depending 
on the topology of $S$, such that the number of subsurfaces with this 
property at each given time is at most $C$. Therefore,
$$d_\T(\sigma, \tau ) \geq \frac 1C \sum |I_Y|.$$ 
Lemmas \ref{lem:length} and \ref{lem:length-ann} imply the desired inequality.
\end{proof}

\subsection{The upper estimate} \label{subsec;upper}
In \cite{minsky:CCII}, 
Masur and Minsky show how to change one marking to another through
elementary moves (described below) efficiently. Their estimate for the 
number of necessary elementary moves closely resembles the estimate in 
\thmref{thm:main}. We use this sequence of elementary moves to 
construct an efficient path connecting two points in $\T(S)$.

There are two types of elementary moves that transform a marking 
$\mu= \{ (\alpha_1, \beta_1), \dots , (\alpha_m,\beta_m) \}$ to a new 
marking.
\begin{enumerate}
\item Twist: Replace $\beta_i$ by $\beta_i'$, where $\beta_i'$ is 
obtained from $\beta_i$ by a Dehn twist or a half twist around $\alpha_i$.
\item Flip: Replace the pair $(\alpha_i,\beta_i)$ with $(\beta_i,\alpha_i)$ 
and, for $j \not = i$, replace $\beta_j$ with a curve $\beta_j'$ that does not 
intersect $\beta_i$, which is now a base curve, in such a way that 
$d_{\alpha_j}(\beta_j, \beta_j')$ is as small as possible
(see \cite{minsky:CCII} for details).
\end{enumerate}
In the first move, a twist can be positive or negative. A half twist
is possible when $\alpha_i$ and $\beta_i$ intersect twice.
The following is a consequence of work done in \cite{minsky:CCII}) and
\cite{minsky:ELCI}.

\begin{proposition} \label{prop:marking-dist}
There exists a large enough $k$ such that:
For markings $\mu$ and $\mu'$, there exists a sequence of markings 
$$\mu=\mu_1, \ldots, \mu_n=\mu',$$
where $\mu_i$ and $\mu_{i+1}$ differ by an elementary move
except, for each $\alpha$ where $d_\alpha(\mu,\mu')\geq k$, there 
is an index $i_\alpha$ so that 
\begin{equation} \label{eq:exponent}
\mu_{i_\alpha+1}= D_\alpha^p \, \mu_{i_\alpha}, 
\quad \text{and} \quad  
|p| \asymp d_\alpha(\mu,\mu').
\end{equation}
Furthermore, 
\begin{equation} \label{eq:n}
n \asymp \sum_{Y\subset S} \big[ d_Y(\mu,\mu') \big]_k,
\end{equation}
where the sum is over all subsurfaces $Y$ that are not annuli.
\end{proposition}

\begin{proof}
We use the definitions and notation used in \cite{minsky:CCII} and
\cite{minsky:ELCI}. \cite[4.6 and 4.20]{minsky:CCII} state that there 
exists a complete hierarchy $H$ whose initial marking is $\mu$ and whose 
terminal marking is $\mu'$.  Any complete marking has a resolution
(\cite[5.4]{minsky:CCII}), that is, there is a sequence of markings
$$\mu = \eta_1, \ldots , \eta_N=\mu'$$
where $\eta_i$ and $\eta_{i+1}$ differ by an elementary move. For $k$ 
large enough, if $d_\alpha(\mu,\mu') \geq k$, the collar of $\alpha$ appears 
as a domain in $H$ (\cite[6.2]{minsky:CCII}) exactly once 
(\cite[5.15]{minsky:ELCI}), and the length of the corresponding geodesic in 
$H$ is comparable to $d_\alpha(\mu,\mu')$ (\cite[6.2]{minsky:CCII}).
That is, the number of twist moves around $\alpha$
used in the resolution is comparable to $d_\alpha(\mu,\mu')$.
The number of the remaining elementary moves is comparable to the sum 
of the lengths of geodesics in $H$ whose domains are not annuli, which is
comparable to (\cite[Lemma 6.2 and Equation (6.4)]{minsky:CCII})
$$ \sum_Y \left[d_Y(\mu, \mu')\right]_k.$$

Our goal is, for any $\alpha$ where $d_\alpha(\mu,\mu')\geq k$,
to rearrange the elementary moves in the resolution so that all the twist 
moves around $\alpha$ are applied consecutively. 
Then we replace the sequence of consecutive twists around $\alpha$ with 
one large step, which is applying $D_\alpha^p$, for some 
$p \asymp d_\alpha(\mu,\mu')$. This will result in the sequence described 
in the statement of the theorem and has the desired length condition. 

We know (\cite[5.16]{minsky:ELCI}) that for every curve $\alpha$, 
the set $J_\alpha$ of indices $i$ such that $\alpha$ is a base curve in 
$\eta_i$ is an interval in $\Z$. Observe that when $\alpha$ is a base curve 
of a marking, a twist move around $\alpha$ and a twist move around any 
other curve can
be rearranged without any complication. The trouble with the flip moves
is that the outcome is not unique. Therefore, after rearranging a flip move
and a twist move, we have to make sure the outcomes of two flip moves 
differ by just a twist around $\alpha$. For example, assume 
$\eta_{i-1}$, $\eta_{i}$ and $\eta_{i+1}$ all contain $\alpha$ as a base curve, 
$\eta_i$ is obtained from $\eta_{i-1}$ by a flip move and 
$\eta_{i+1} = D_\alpha \eta_i$. Then, replacing $\eta_i$ with
$\eta_i' = D_\alpha \, \eta_{i-1}$ in our sequence will result in
a sequence that is still a resolution of $H$. Because $\eta_i$
is obtained from $\eta_{i-1}$ by applying a flip move, 
$D_\alpha \, \eta_i$ is also obtained from $D_\alpha \, \eta_{i-1}$ by a flip 
move ($D_\alpha$ is a homeomorphism).
Therefore, we can rearrange the elementary moves in $J_\alpha$
so that all the twist moves around $\alpha$ are done consecutively.
\end{proof}

\begin{remark} \label{rmk:k}
The constant $k$ can be chosen as large as necessary, and
the constants involved in \eqref{eq:n} depend on $k$ and the 
topology of $S$ (see Theorem 6.12 in \cite{minsky:CCI}).
Therefore, we can assume $k \geq k_0$, where $k_0$
is as chosen in \propref{prop:lower}.
\end{remark}

For a marking $\mu$, let $\short(\mu)$ be the set of points
in $\T(S)$ where all curves in $\mu$ have hyperbolic length
less than $L_0$ ($L_0$ as on page \pageref{page:L}).  This is a compact subset of $\T(S)$. We define $f(\mu,\mu')$ to be the maximum distance between 
an element in $\short(\mu)$ and an element in $\short(\mu')$.
 
\begin{lemma} \label{lem:step}
If $i=i_\alpha$, where $\alpha$ is a curve with $d_\alpha(\mu,\mu')\geq k$, then
$$f(\mu_i, \mu_{i+1}) \asymp \log d_\alpha(\mu,\mu').$$
Otherwise,
$$ f(\mu_i,\mu_{i+1}) = O(1).$$
\end{lemma}

\begin{proof}
Since $\short(\mu)$ is compact, it is enough to bound the minimum
distance between $\short(\mu_i)$ and $\short(\mu_{i+1})$.

Assume $i=i_\alpha$, for $\alpha$ as above, and let $\sigma$ be a point
in $\short(\mu_i)$. Then, for some $|p| \asymp d_\alpha(\mu,\eta)$,
$\tau= D^p_\alpha \, \sigma$ is a point in $\short(\mu_{i+1})$.
The lengths of $\alpha$ in $\sigma$ and $\tau$ are less than $L_0$, 
therefore, $\sigma$ and $\tau$ are bounded distance from 
points $\sigma'$ and $\tau'=D_\alpha^p(\sigma')$, where the lengths of 
$\alpha$ in $\sigma'$ and $\tau'$ are less than $\ep_0$. Taking $\Gamma=\{\alpha\}$ and $\pi$ as 
in \thmref{thm:product}, the following holds: the distance between $\sigma'$ and 
$\tau'$ equals, up to additive error, the distance in $\hyp_\alpha$ between
$\pi_\alpha(\sigma')$ and $\pi_\alpha(\tau')$, which, up to 
multiplicative error, equals $\log |p|$. Therefore, the distance between
$\sigma$ and $\tau$ is comparable to $\log |p|$. 

Otherwise, $\mu_i$ and $\mu_{i+1}$ differ by an elementary move.
Note that there are only finitely many such pairs of markings up to homeomorphism. Therefore, there exists a uniform upper bound for the 
minimum distance between  $\short(\mu_i)$ and $\short(\mu_{i+1})$,
depending on the topology of $S$ only.
\end{proof}

\begin{proposition} \label{prop:upper}
Let $\sigma_1$, $\sigma_2$ be in the $\ep_o$--thick part of $\T(S)$ and $\mu_1$ and $\mu_2$ be the short-markings in $\sigma_1$ and 
$\sigma_2$, respectively. Then
\begin{equation*}
d_\T(\sigma_1, \sigma_2 ) \prec \sum_Y \big[ d_Y(\mu_1, \mu_2) \big]_k +
\sum_\alpha \log \big[ d_\alpha(\mu_1,\mu_2) \big]_k.
\end{equation*}
\end{proposition}
\begin{proof}
Let $\mu_1=\bar \mu_1, \dots, \bar \mu_n=\mu_2$ be the path in $\M(S)$ 
described in \propref{prop:marking-dist}. For each $i$, let $\sigma_i$ be a
point in $\short(\bar \mu_i)$ and let $g_i$ be the geodesic segment 
connecting $\sigma_i$ to $\sigma_{i+1}$. The distance in $\T(S)$ between
$\sigma_1$ and $\sigma_2$ is less than the sum of the lengths of the $g_i$.
\lemref{lem:step} states that the lengths of the $g_i$ are uniformly bounded
except when $i=i_\alpha$ and $d_\alpha(\mu_1,\mu_2) \geq k$, in
which case the length of $g_i$ is comparable with 
$\log d_\alpha(\mu_1,\mu_2)$.
Therefore,
$$d(\sigma,\tau) \prec n \, O(1) + 
\sum_\alpha \log  \big[d_\alpha(\mu_1,\mu_2) \big]_k.$$
\propref{prop:marking-dist} finishes the proof.
\end{proof}

\begin{proof}[Proof of \thmref{thm:main}]
Propositions \ref{prop:lower} and \ref{prop:upper} provide a lower estimate
and an upper estimate for the distance between $\sigma$ and $\tau$.
Since $k \geq k_0$ (see \remref{rmk:k}), the estimate given in 
\propref{prop:upper} is smaller than the one given in \propref{prop:lower}.
Therefore $d_\T(\sigma_1, \sigma_2)$ is comparable to
\begin{equation*}
\sum_Y \big[ d_Y(\mu_1, \mu_2) \big]_k +
\sum_\alpha \log \big[ d_\alpha(\mu_1,\mu_2) \big]_k. \qedhere
\end{equation*}
\end{proof}

\section{The general case\label{sec:general}}
\noindent
In this section we give an estimate for the distance between two arbitrary 
points in the \Teich space. Let  $\sigma_1$ and $\sigma_2$ be two points in $\T(S)$ and $g \co [a,b] \to \T(S)$ be the geodesic arc connecting them. 
If $\sigma_1$ and $\sigma_2$ are not in the thick part of $\T(S)$, then the 
set of short curves 
in $\sigma_1$ and $\sigma_2$ does not contain enough information to allow
us to 
estimate the distance between $\sigma_1$ and $\sigma_2$; we also need to know how short these curves are.
Therefore, our estimate for the distance contains terms measuring
the distance between $\sigma_1$ and $\sigma_2$ and the thick part of \Teich 
space. An additional complication arises from the case where a curve is short in 
both $\sigma_1$ and $\sigma_2$ and remains short along the geodesic.
However, the basic idea behind both \thmref{thm:main} and \thmref{thm:general} 
is that efficient paths in the space of markings are closely related
to geodesics in \Teich space.

Let $\ep_0$ be as before. Define $\Gamma$ to be the set of curves that are short in 
both $\sigma_1$ and $\sigma_2$, and, for $i=1,2$, define $\Gamma_i$ to be the 
set of curves that are short in $\sigma_i$ but not in $\sigma_{3-i}$.
Let $\mu_1$ and $\mu_2$ be short-markings on $\sigma_1$ 
and $\sigma_2$, respectively.

\begin{theorem} \label{thm:general}
The distance in $T(S)$ between $\sigma_1$ and $\sigma_2$ is given
by the following formula:
\begin{equation} \label{eq:big} \begin{split}
d_\T(\sigma_1, \sigma_2 ) \asymp 
& \sum_Y \big[ d_Y(\mu_1, \mu_2) \big]_k +
\sum_{\alpha \not \in \Gamma} \log \big[ d_\alpha(\mu_1,\mu_2) \big]_k + \\
& + \max_{\alpha \in \Gamma} \, d_{\hyp_\alpha}(\sigma_1,\sigma_2) +
\max_{\begin{subarray}{l} \alpha \in \Gamma_i \\ i=1,2 \end{subarray}}
\, \log \frac 1{l_{\sigma_i}(\alpha)}.
\end{split} 
\end{equation}
\end{theorem}

\begin{proof}
\thmref{thm:product} implies that
$$d_\T(\sigma_1,\sigma_2) \asymp 
\max_{\alpha \in \Gamma} \, d_{\hyp_\alpha}(\sigma_1,\sigma_2)
+ d_{\T(S \setminus \Gamma)} \big( \pi_0(\sigma_1),\pi_0(\sigma_2) \big).$$
This accounts for the third term on the right-hand side of 
Equation \eqref{eq:big}.  Therefore, without loss of generality,
we can assume $\Gamma = \emptyset$.

Let $\sigma_1'$ and $\sigma_2'$ be points in the thick part of the \Teich space
that have the same short-markings as $\sigma_1$ and $\sigma_2$. 
We have:
$$d_\T(\sigma_1, \sigma_2) \leq d_\T(\sigma_1,\sigma_1') +
 d_\T(\sigma_1',\sigma_2') +d_\T(\sigma_2',\sigma_2),$$
The sum of the first two terms in \eqref{eq:big} is comparable with
$d_\T(\sigma_1',\sigma_2')$. Also,
$$ 
d_\T(\sigma_1,\sigma_1') \asymp 
\max_{\beta \in \Gamma_1} \, \log \frac 1{l_{\sigma_1}(\beta)}
\quad \text{and} \quad 
d_\T(\sigma_2,\sigma_2') \asymp 
\max_{\gamma \in \Gamma_2} \, \log \frac 1{l_{\sigma_2}(\gamma)}.
$$
Therefore, the right side of \eqref{eq:big} is an upper bound for 
$d_\T(\sigma_1,\sigma_2)$ (up to additive and multiplicative constants). 

To show that the right side of \eqref{eq:big} is also a lower bound
for $d_\T(\sigma_1,\sigma_2)$, we follow the same argument as in 
\secref{subsec:lower}.  However, we can not use \lemref{lem:length-ann}
when $\alpha$ is short in either $\sigma_1$ or $\sigma_2$ and
using the previous argument we can conclude only that
\begin{equation} \label{eq:almost}
d_\T(\sigma_1, \sigma_2 ) \succ \sum_Y \big[ d_Y(\mu_1, \mu_2) \big]_k +
\sum_{\alpha \not \in \Gamma\cup\Gamma_1\cup \Gamma_2} 
\log \big[ d_\alpha(\mu_1,\mu_2) \big]_k.
\end{equation}

For every $\alpha \in \Gamma_1$, we have
$$d_\T(\sigma_1,\sigma_2) = \log K(\sigma_1,\sigma_2) \geq
\log \left| \frac {\Ext_{\sigma_2}(\alpha)}{\Ext_{\sigma_1}(\alpha)} \right|
\succ \log \frac 1{l_{\sigma_1}(\alpha)}.$$
A similar statement is true for $\alpha \in \Gamma_2$. Hence
\begin{equation} \label{eq:fourth-term}
d_\T(\sigma_1, \sigma_2 ) \succ  
\max_{\begin{subarray}{l} \alpha \in \Gamma_i \\ i=1,2 \end{subarray}}
\, \log \frac 1{l_{\sigma_i}(\alpha)}.
\end{equation}
It remains to show, for $\alpha \in \Gamma_1\cup \Gamma_2$, 
that $d_\T(\sigma_1,\sigma_2) \succ \log d_\alpha(\mu_1,\mu_2)$.
Let $\beta_1$ and $\beta_2$ be the transverse curves to $\alpha$
in $\mu_1$ and $\mu_2$. We know
$$|tw_{\sigma_1}^+ - tw_{\sigma_2}^+|
= |i_\alpha(\np,\beta_1)-i_\alpha(\np,\beta_2)|
\eadd i_\alpha(\beta_1,\beta_2)=d_\alpha(\mu_1,\mu_2).$$
Therefore, it is sufficient to show that 
$d_\T(\sigma_1,\sigma_2) \succ \log |tw_{\sigma_1}^+ - tw_{\sigma_2}^+| $.
\thmref{thm:twist-comparison} implies that
$$
 |tw_{q_1}^+ - tw_{q_2}^+| 
 \succ \left| tw_{\sigma_1}^+ - tw_{\sigma_2}^+ - 
  O(\frac 1{l_{\sigma_1}(\alpha)} + \frac 1{l_{\sigma_2}(\alpha)}) \right|;
$$
therefore,
$$|tw_{q_1}^+ - tw_{q_2}^+| +  \frac 1{l_{\sigma_1}(\alpha)} + \frac 1{l_{\sigma_2}(\alpha)} \succ |tw_{\sigma_1}^+ - tw_{\sigma_2}^+|,$$
and Equation \eqref{eq:q-twisting} implies that the $q_t$--twisting
parameter changes at most exponentially fast; hence,
$$d_\T(\sigma_1,\sigma_2) \succ \log (tw_{q_1}^+ - tw_{q_2}^+).$$  
We also know that 
$$d_\T(\sigma_1,\sigma_2) \succ \log \frac 1{l_{\sigma_1}(\alpha)}
\quad\text{and}\quad
d_\T(\sigma_1,\sigma_2) \succ \log \frac 1{l_{\sigma_2}(\alpha)}.$$
From the last three equations, we can conclude
$$
d_\T(\sigma_1,\sigma_2)  \succ \log | tw_{\sigma_1}^+ - tw_{\sigma_2}^+ | 
.$$
Therefore,
\begin{equation} \label{eq:second-term}
\forall \alpha \in \Gamma_1 \cup \Gamma_2, \qquad
d_\T(\sigma_1,\sigma_2)  \succ \log d_\alpha(\mu_1,\mu_2). 
\end{equation}
The combination of Equations \eqref{eq:almost}, \eqref{eq:fourth-term}
and \eqref{eq:second-term} provides the desired lower bound and finishes
the proof.
\end{proof}

\bibliographystyle{alpha}
\bibliography{../main}

\end{document}